\documentclass[11pt]{amsart}
\usepackage{amssymb,latexsym,amscd, amsthm, epsfig}
\usepackage{color}
\newtheorem{Thm}{Theorem}[section]

\newtheorem{Prop}[Thm]{Proposition}
\newtheorem{Lem}[Thm]{Lemma}
\theoremstyle{definition}
\newtheorem{Def}[Thm]{Definition}
\newtheorem{Ex}[Thm]{Example}
\newtheorem{Rmk}[Thm]{Remark}

\date{}

\newcommand{\PP}{{\mathbb{P}}}
\newcommand{\CC}{{\mathbb{C}}}

\newcommand{\triplet}[3]{\PP^{#1}\times\PP^{#2}\times\PP^{#3}}

\begin{document}

\title[On generic identifiability of $3$-tensors of small rank]{On generic identifiability of $3$-tensors\\ of small rank}
\author{Luca Chiantini}
      \address{Luca Chiantini\\
Universit\'a degli Studi di Siena\\
     Dipartimento di Scienze Matematiche e Informatiche\\
     Pian dei Mantellini, 44\\
     I -- 53100 Siena
     }
   \email{chiantini@unisi.it}
  \urladdr{http://www.mat.unisi.it/newsito/docente.php?id=4}
\author{Giorgio Ottaviani}
      \address{Giorgio Ottaviani\\
Universit\'a degli Studi di Firenze\\
     Dipartimento di Matematica ''U. Dini''\\
     Viale Morgagni, 67a\\
     I -- 50134 Firenze
     }
   \email{ottavian@math.unifi.it}
  \urladdr{http://web.math.unifi.it/users/ottavian/}
  
\begin{abstract}
We introduce an inductive method for the study of the 
uniqueness of  decompositions of
tensors, by means of tensors of rank $1$. The method is based on the geometric notion
of {\it weak defectivity}. For  three-dimensional tensors
of type $(a,b,c)$, $a\leq b\leq c$, our method proves that
the decomposition is unique (i.e. $k$-identifiability holds) for general tensors 
of rank $k$, as soon as $k\leq (a+1)(b+1)/16$. This improves considerably the known range 
for identifiability. The method applies also to tensor of higher dimension. 
For tensors of small size, we give a complete list of
situations where identifiability does not hold. Among them,
there are $4\times 4\times 4$ tensors of rank $6$, an interesting case
because of its connection with the study of DNA strings.
\end{abstract}
\maketitle

\section{Introduction}
\subsection{Statement of main results}

Let $A$, $B$, $C$ three complex vector spaces, 
of dimension $a$, $b$, $c$ respectively. 
A tensor $t\in A\otimes B\otimes C$ is said to have {\it rank} $k$ 
if there is a decomposition
$$t=\sum_{i=1}^k u_i\otimes v_i\otimes w_i$$ 
with $u_i\in A, v_i\in B, w_i\in C$ and the number of summands
$k$ is minimal. Such a decomposition is said to be {\it unique} 
if for any other expression
$$t=\sum_{i=1}^k u'_i\otimes v'_i\otimes w'_i$$
there is a permutation $\sigma$ of $\{1,\ldots , r\}$  such that 
$$u_i\otimes v_i\otimes w_i=u'_{\sigma(i)}\otimes v'_{\sigma(i)}\otimes w'_{\sigma(i)}\quad\forall i=1,\ldots, k.$$ 
When $t$ has a unique  decomposition, the vectors $u_i\in A$, 
$v_i\in B$, $w_i\in C$ can be {\it identified} uniquely from $t$, 
up to scalars. 

It is known that the set of tensors of rank $k$
consists of a dense subset of an irreducible 
algebraic variety $S_k(Y)$,  which is called the {\it $k$-th secant variety}
of the variety $Y$ of tensors of rank one. 
This last variety is isomorphic to the (cone over the) 
Segre product $\PP(A)\times\PP(B)\times\PP(C)$.

The main result of our paper determines a bound for the rank, 
in terms of the dimensions of the
vector spaces, which implies identifiability

\begin{Thm}\label{main} 
 Let $a\le b\le c$. Let $\alpha$, $\beta$ be maximal such 
 that $2^\alpha\le a$  and $2^\beta\le b$.
The general tensor  $t\in A\otimes B\otimes C$ of rank $k$ has a 
unique decomposition, if $k\le 2^{\alpha+\beta-2}$.
\end{Thm}

So if $a$, $b$ are both a power of $2$, then the general tensor of 
rank $k$ has a unique decomposition if $k\le\frac{ab}{4}$. 
In the general case, the inequality of the theorem can be written as
$k\le 2^{\left(\lfloor \log_2{a}\rfloor+\lfloor \log_2{b}\rfloor-2\right)}$.
Since $\frac{a+1}{2}\le 2^\alpha$ and 
$\frac{b+1}{2}\le 2^\beta$, one can say that the unique 
decomposition holds if $k\le (a+1)(b+1)/16$.

In our terminology, when the unique decomposition holds
for the general tensor of rank $k$, we will say that
the variety of tensors of rank one is {\it $k$-identifiable}. 

Here the meaning of ``general'' is that, among tensors of rank $k$,
the ones which do not have a unique decomposition 
consist in a set of zero measure, more specifically
in a proper subvariety of $S_k$.

In particular, the Theorem applies to ``cubic'' tensors. 
The general tensor  $t\in A\otimes A\otimes A$ of rank 
$k$ has a unique decomposition if $k\le \frac{a^2}{16}$
(indeed, the Theorem provides a better bound, when
$a$ is close to a power of $2$).

Our bound is log-asymptotically sharp, in the following sense.
As explained in Proposition \ref{maxk}, one cannot
have a unique decomposition, when the rank exceeds a value $k_{max}=
k(a,b,c)$, which depends on $a,b,c$. 
Then $\sup_c \frac{k(a,b,c)}{ab}$ is finite.
On the other hand, even for tensors of small size, 
the result is not sharp. In the first cases, with the help of a computer, 
we can  improve Theorem \ref{main}. 

Unique decomposition has been studied by several authors, and there
is a huge amount of literature, on this theme.
Let us remind that Strassen and Lickteig (\cite{Lick}) proved that
the general tensor $t\in A\otimes A\otimes A$ has rank 
$\lceil \frac{a^3}{3a-2}\rceil$ for $a\neq 3$
and rank $5$ for $a=3$  (indeed, the case $a=3$ is known to be 
{\it defective}, meaning that the corresponding $4$-secant variety
has dimension smaller than the expected one). 
In this case, the aforementioned bound implies that, if
$a\geq 3$, then the generic 
tensor of rank $k$ can have a unique decomposition only
if $k\le \lceil \frac{a^3}{3a-2}\rceil -1$.
The following theorem shows that this bound is almost always 
achieved, for small $a$.

\begin{Thm}\label{main2}
The general tensor  $t\in A\otimes A\otimes A$ of rank $k$ 
has a unique decomposition if $k\le k(a)$ where
$$\begin{array}{cc|rrrrrr|rrr}&a&2&3&4&5&6&7&8&9&10\\
\hline\\
&k(a)&2&3&5&9&13&18&22&27&32
\end{array}$$
\end{Thm}

\vskip.5cm
A more general list, which holds in the non cubic case, is given 
in section \ref{fin}.

Comparing the previous table with the table 
of the general rank (for $a> 3$, the general rank $-1$ 
is the best possible achievement),
and with Kruskal's result (see Proposition \ref{Kruskal}), 
one can appreciate the improvement.

$$\begin{array}
{cc|rrrrrr|rrr}&a&2&3&4&5&6&7&8&9&10\\
\hline\\
\textrm{gen.rank\ } (a\neq 3)\ \textrm{\cite{Lick}}&\lceil \frac{a^3}{3a-2}\rceil&2&4&7&10&14&19&24&30&36\\
\hline\\
\textrm{Kruskal bound\ \cite{Kru}}&\lfloor\frac{3a-2}{2}\rfloor&2&3&5&6&8&9&11&12&14
\end{array}$$

\vskip.5cm
The more evident lack of uniqueness is when $a=4$ and $k=6$.
The case $a=4$ is particularly interesting due to the 
models in phylogenetics \cite{AR,ERSS},
where a basis in $\CC^4$ can be indexed by the nucleotids 
$\{A, C, G, T\}$.

\begin{Thm}\label{main3}
The general tensor $t\in\CC^4\otimes\CC^4\otimes\CC^4$ of rank 
$6$ has exactly two decompositions. 
\end{Thm}

It is interesting that the exception on uniqueness ($a=4$) holds 
very close to the defective case $a=3$. This phenomenon is quite 
general and it can be already encountered in the case of symmetric tensors.

\subsection{A few historical remarks}
In this subsection we sketch how our result fills in the literature.

The most celebrated result about uniqueness of decomposition of 
tensors is due to Kruskal \cite{Kru}. It is often quoted in 
terms of Kruskal's rank. A consequence of Kruskal's criterion 
is the following statement, which applies to generic tensors (see 
 Corollary 3 in \cite{AMR}).

\begin{Prop}\label{Kruskal} ({\bf Kruskal's criterion})
The generic tensor $t\in A\otimes B\otimes C$ of rank $k$
 has a unique decomposition if 
$$ k \leq \frac 12 \left[\min(a,k)+\min(b,k)+\min(c,k)-2\right]$$
In the cubic case,  the generic tensor $t\in A\otimes A\otimes A$ 
of rank $k$ has a unique decomposition if 
 $$ k\le \frac{3a-2}{2}$$
\end{Prop}

Kruskal's result is so important in the literature,
that recently there have been published (at least!) three different proofs
\cite{Land,R,St}.

De Lathauwer (\cite{Lat}) proves that the generic tensor 
$t\in A\otimes B\otimes C$ of rank $k$ has a unique decomposition if 
 $ k\le c$ and $k(k-1)\le a(a-1)b(b-1)/2$.
Rhodes, in \cite{R} addresses explicitly, as a problem at the end of the 
introduction, the need of sufficient conditions,
stronger than Kruskal's, that guarantee the uniqueness of 
the decomposition, for generic tensors. Our Theorem \ref{main} gives a sufficient 
condition which improves both Kruskal's and de Lauthawer's bounds.

The tensor decomposition we are looking for are called also 
Candecomp or Parafac decompositions in the numerical literature.
Among recent surveys on the topic, see \S 3.2 in \cite{KB} and 
Landsberg book \cite{Land0}, which tries to use a language 
understandable by both the numerical and the geometrical communities.
From this point of view, one should also consider section 2 of
\cite{AMR}, an interesting bridge between the two worlds.

\subsection{Outline of the proof}
In a line, our technique consists in putting together the 
inductive approach of \cite{AOP} with the tool of weak defectivity 
developed in \cite{CC1} and \cite{CC2}.

We consider the projective space of tensors $\PP(A\otimes B\otimes C)$.
In this space, the tensors of rank one give the Segre variety 
$\PP(A)\times\PP(B)\times\PP(C)$. 

Our geometric point of view consists in the use of the celebrated 
Terracini's Lemma, which allows to study the identifiability of varieties,
using properties of their tangent spaces. We refer to \cite{CC1}
and \cite{CC2} for a more precise account of the theory behind. 

A variety is called {\it tangentially $k$ weakly defective}
($k$-twd, see Definition \ref{twd}) if  the span of the tangent
spaces at $k$ general points of $X$, is tangent also in some 
other points.

It is a consequence of Terracini's Lemma that, if $X$ is $k$-not 
twd, then the general tensor of rank $k$  has a unique decomposition.

So our aim is to prove the $k$-not twd of Segre varieties 
$X=\PP(A)\times\PP(B)\times\PP(C)$.
The proof is performed by induction, by splitting
$A=A'\oplus A''$ and by specializing some points on the lower 
dimensional Segre varieties $\PP(A')\times\PP(B)\times\PP(C)$ and
$\PP(A'')\times\PP(B)\times\PP(C)$. It turns out that the 
induction works if we prove a stronger statement, concerning 
the so called $(k,p,q,r)$-weakly defectivity,
which is defined in section \ref{indsec}.

\subsection{Outline of the paper}
In section \ref{prel} we develop the basic notations on Segre varieties and 
weak defectivity. At the end of this section we prove the cases $a\le 7$
of the Theorem \ref{main2}. Section \ref{indsec} contains  the definition \ref{kpqr} 
of $(k,p,q,r)$-defectivity and the inductive step (Prop. \ref{induz}). 
At the end of this section we prove the remaining cases 
of the Theorem \ref{main2}.
In the section \ref{pfs} we prove the Theorem \ref{main}.
In section \ref{fin} we prove the Theorem \ref{main3} and we give other 
examples of small dimension. Also we expose a list of all the examples
of triple Segre product that we know when the uniqueness 
for general tensors of a given rank does not hold. 
In section \ref{many}, we show an extension of the
previous results to products of many factors.

\section{Preliminaries on Segre varieties}\label{prel}

Let $A,B,C$ be complex vector spaces, of dimension $a,b,c$ respectively.
Consider the product  $X=\PP(A)\times\PP(B)\times\PP(C)$. 
$X$ is naturally embedded, by means of the Segre map, into
$\PP^N$, where $N=abc-1$.

Sometimes, when there is no need to specify the vector spaces, 
we will refer to the variety $X$ also as $\triplet {a-1}{b-1}{c-1}$. 

Call $S^k(X)$ the $k$-th secant variety of $X$, defined as the closure
of the union of linear spans of $k$ general points in $X$.

\begin{Def} 
$X$ is called $k$-identifiable if a \emph{general} element
in $S_k(X)$ has a unique expression as sum of $k$ elements in $X$.

From the tensorial point of view, this means that a general tensor
of type $a\times b\times c$ and rank $k$, can be written uniquely 
(up to scalar multiplication) as a sum of $k$ decomposable tensors.
\end{Def}

\begin{Prop}\label{maxk} There is a maximal rank for which 
the $k$-identifiability of tensors is possible, namely
$$k_{max}=  \lfloor \frac {N+1}{\dim(X)+1} \rfloor= \lfloor  
\frac{abc}{a+b+c-2} \rfloor.$$
\end{Prop}
\begin{proof}
For $k>k_{max}$, the abstract secant variety
$$ AbS^k(X)=\{(x_1,\dots,x_k,u)\in X^k\times\PP^N: 
u\in<x_1,\dots,x_k>\}$$ 
has dimension bigger than $N$, so that necessarily the general 
$u\in S^k(X)$ belongs to infinitely many $k$-secant spaces.
\end{proof}

Our theoretical starting point is a criterion for $k$-identifiability,
which follows from the Terracini's Lemma, 
which we will use under the following form (see e.g. \cite{CC1})

\begin{Lem} ({\bf Terracini}) Let $X$ be an irreducible variety and 
consider a general point $u\in S_k(X)$. If $u$ belongs to the span
of points $x_1,\dots, x_k\in X$, then the tangent space to $S_k$ at $u$
is the span of the tangent spaces to $X$ at the points $x_1,\dots,x_k$.
\end{Lem}

Our criterion is the following: 
 
\begin{Prop}\label{criterio}
Let $X\subset \PP^N$ be a non-degenerate, irreducible variety 
of dimension $n$. Consider the following statements:

(i) $X$ is $k$-identifiable

(ii) Given $k$ general points $x_1,\ldots x_k\in X$, then the 
span $<T_{x_1}X,\ldots, T_{x_k}X>$
contains $T_xX$ only if $x=x_i$ for some $i=1,\ldots k$.

(iii)  there exists a set of $k$ particular points 
$x_1,\ldots x_k\in X$, such  that the span 
$<T_{x_1}X,\ldots, T_{x_k}X>$ contains $T_xX$ only if $x=x_i$ 
for some $i=1,\ldots k$.

Then we have (iii) $\Longrightarrow$ (ii) $\Longrightarrow$ (i).
\end{Prop}
\begin{proof} (iii) $\Longrightarrow$ (ii) follows 
at once by semicontinuity.

Let us prove that (ii) $\Longrightarrow$ (i). 
Take a general point $u\in S_k(X)$ and assume that $u$ belongs
to  the span of points $x_1,\dots, x_k\in X$. By the generality of $u$, 
we may assume that $x_1,\dots, x_k$ are general points of $X$.
If $u$ also belongs to the span of points $y_1,\dots, y_k\in X$,
with at least one of them, say $y_1$, not among the $x_i$'s,
then, by Terracini's Lemma, the span of the tangent spaces to $X$ at
the points $x_i$'s, which is the tangent space to $S_k(X)$ at $u$, also
contains the tangent space to $X$ at $y_1$. This contradicts (ii). 
\end{proof}

Condition (ii) of the previous Proposition is strongly
related with the notion of $k$-weak defectivity.

In \cite{CC1}, C. Ciliberto and the first author give the following
definition: a variety $X$ is $k$-weakly defective if the general
\emph{hyperplane} which is tangent to $X$ at $k$ general points 
$x_1,\dots,x_k$, is also tangent in some other point $y\neq
x_1,\dots, x_k$.

It is clear that a variety which does not satisfy condition (ii)
of the Proposition, is also $k$-weakly defective. However the
converse does not hold.

\begin{Ex} \label{nonbitangent}Consider the Segre product $X=\PP^1\times\PP^2$.
It is classical (see e.g. Zak's Theorem on tangencies in \cite{zak})
that the tangent space at one point to a smooth variety is not
tangent elsewhere. 

On the other hand, a general hyperplane tangent to $X$ at one point, 
is also tangent along a line. Indeed, it is well known that
the dual variety of $X$ is not a hypersurface (see \cite{Ein}).
Thus $X$ is $1$-weakly defective. 
\end{Ex}

For maintaining the consistency with all the previous
notation in this subject, we dare proposing the following:

\begin{Def}\label{twd} If $X$ satisfies condition (ii) of the previous
Proposition, we will say that $X$ is \emph{ $k$-not tangentially 
weakly defective}. Otherwise, we say that $X$ is \emph{$k$-tangentially weakly defective} 
($k$-twd, for short).
\end{Def}
 
We understand that the notation is becoming odd. However,
the increasing number of definitions is a phenomenon which
also occurs in the study of \emph{contact loci}, which
seems however helpful for applications to the Geometry
of secant varieties (see e.g. \cite{CC3}).
 
Weak defectivity has been intensively studied in \cite{CC1}.
Notice that when $X$ is $k$ weakly defective, then a general hyperplane
tangent to $X$ at general points $x_1,\dots,x_k$ is also
tangent along a positive dimensional variety. We do
not know if a similar phenomenon takes place also 
for $k$-twd.
 
Relations between $k$-weak defectivity and $k$-twd 
are probably stronger than expected, at least as far as
one is interested in $k$-identifiability. We do not develop
further this analysis.

Notice than, when we deal with inductive steps in the proofs,
we will need an even more complicated notion of
weak defectivity. Compare with Definition \ref{kpqr} below.
\smallskip

For our purposes, Proposition \ref{criterio} establishes that
 $k$-not tangentially weakly defectivity implies
$k$-identifiability, 
when $N\geq k(n+1)$.

\begin{Rmk}
Let us notice that, by Proposition \ref{maxk}, 
if $N+1<(k+1)n$, then $k$-identifiability is excluded.
Thus, the criterion of Proposition \ref{criterio} cannot be applied  only for
at most one value of $k$, namely $k=(N+1)/(n+1)$,
which occurs only when $N+1$ is a multiple of $n+1$.
E.g., our criterion could not be applied to
study the $2$-identifiability of $\triplet 111$.
\end{Rmk} 

Now we are already able to prove
the first cases of Theorem \ref{main2}.

{\it Proof of the Theorem \ref{main2} in case $a\le 7$.}

The proof is a straightforward application of Proposition \ref{criterio}. 
A random choice of $k(a)$ points satisfies condition (iii) of
Proposition \ref{criterio}. Then $X$ is $k$-identifiable.
The Macaulay2 files which we used are available as ancillary files in the arXiv submission of this paper.

\begin{Rmk} More powerful computers and/or better suited algorithms will allow eventually to check the condition (iii) for larger values of $a$,
and we encourage experts in Numerical Algebraic Geometry in going further.
We stopped at $a=7$, because for $a=8$
our algorithm on a common PC consumed too much time and memory. In the next section we
show how the computation for larger values of $a$ can be reduced to other computations for smaller values of $a$.
\end{Rmk}

\section{The inductive statement}\label{indsec}

The inductive criterion makes use of the fact that
if $x=u\otimes v\otimes w$ is a point of 
$X=\PP(A)\times\PP(B)\times\PP(C)$, 
then the tangent space $T_xX$ is the projectification
of the linear space $A\otimes v\otimes w +
u\otimes B\otimes w + u\otimes v\otimes C$.

The idea is to fix two linear subspaces $A',A''$ of $A$,
such that $A=A'\oplus A''$, then split the set of $k$ points in two
subsets and specialize them to the two spaces
$\PP(A')\times\PP(B)\times\PP(C)$ and $\PP(A'')\times\PP(B)\times\PP(C)$.
Then, the implication (iii) $\Longrightarrow$ (i) of 
Proposition \ref{criterio} suggests that one could play induction.

Unfortunately, the situation is a little bit more complicated,
since one cannot translate condition (ii) of Proposition \ref{criterio}
into the analogous condition on lower-dimensional spaces. 

Instead, following the idea of \cite{AOP} (Theorem 3.4)
(suggested also from the Splitting Method of \cite{BCS}), 
we need a more elaborated condition.

\begin{Def}\label{kpqr} A triple product 
$X=\PP(A)\times\PP(B)\times\PP(C)$ 
is called {\it $(k,p,q,r)$-not weakly defective} if:
 
for $k$ general points $x_1,\ldots x_k\in X$,

for $p$ general points $u_i\in\PP(B)\times\PP(C)$,

for $q$ general points $v_i\in\PP(A)\times\PP(C)$,

for $r$ general points $w_i\in\PP(A)\times\PP(B)$,

then the span of 
$T_{x_i}X$, $A\otimes u_i$, $B\otimes v_i$, $C\otimes w_i$
contains  $T_xX$ if and only if $x=x_i$, for some $i=1,\ldots k$.
Otherwise $X$ is called {\it $(k,p,q,r)$-weakly defective}.

Clearly, $(k,0,0,0)$ weak defectivity
coincides with $k$-twd.
\end{Def}

\begin{Rmk}
We will often use the computer algorithm,
available in our arXiv submission,
to prove that some triple Segre product is
$(k,p,q,r)$-not weakly defective.
 
For instance, the algorithm shows that 
$\triplet 222$ is $(1,2,1,1)$-not and $(2,1,1,1)$-not weakly defective.
This is rather interesting, because
$\PP^2\times\PP^2\times\PP^2$ is $3$-defective.
\end{Rmk}

\begin{Ex} Consider $A$, $B$, $C$, all of dimension $2$ with basis
$\{u_1,u_2\}$, $\{v_1,v_2\}$, $\{w_1, w_2\}$.
Then $T_{u_1v_1w_1}+u_2v_2C=T_{u_2v_2w_1}+u_1v_1C$.
This shows that $\PP^1\times\PP^1\times\PP^1$ is $(1,0,0,1)$ 
weakly defective. Nevertheless, $T_{u_1v_1w_1}+u_2v_2C$
has the expected (affine dimension) $6$ and it does not fill the ambient space.
\end{Ex}

\begin{Rmk}\label{ovv} (a) With the previous notation, by semicontinuity
it is clear that when $X$ is  $(k,p,q,r)$-not weakly defective,
then it is also  $(k,p,q,r)$-not weakly defective,
whenever $(k',p',q',r')\leq (k,p,q,r)$, in the strict ordering.

(b) By semicontinuity, $X$ is  $(k,p,q,r)$-not weakly defective
whenever one gets that for {\it particular} sets of points
$\{x_i\}$, $\{u_i\}$, $\{v_i\}$ and $\{w_i\}$ as above,
then the span of 
$T_{x_i}X$, $A\otimes u_i$, $B\otimes v_i$, $C\otimes w_i$
contains  $T_xX$ if only if $x=x_i$, for some $i=1,\ldots k$.

(c) By Proposition \ref{criterio}, one gets soon that
 $(k,0,0,0)$-not weakly defective implies  $k$-identifiable.
\end{Rmk}

We will often apply the following reduction step:

\begin{Lem} \label{riduz} 
Assume that $\triplet {a-1}{b-1}{c-1}$ is $(k,p,q,r)$-not
weakly defective. Then $\triplet{a'}{b'}{c'}$ is $(k,p,q,r)$-not
weakly defective, for any triple $(a',b',c')>(a-1,b-1,c-1)$
(in the strict ordering).
\end{Lem}
\begin{proof} We need just to prove the statement for
$(a',b',c')=(a,b-1,c-1)$. Write 
$X=\triplet a{b-1}{c-1}=\PP(A')\otimes\PP(B)\otimes\PP(C)$
so that $\dim(A')=a+1$.

Assume that $X$ is $(k,p,q,r)$-weakly defective.
Thus, for $k$ general points $x_1,\dots, x_k\in X$, 
$p$ general points $u_i\in\PP(B)\times\PP(C)$,
$q$ general points $v_i\in\PP(A')\times\PP(C)$,
$r$ general points $w_i\in\PP(A')\times\PP(B)$,
then the span $\Lambda$ of the tangent spaces to $X$ at the $x_i$'s
and the spaces $A'\otimes u_i$, $B\otimes v_i$, 
$C\otimes w_i$, is also tangent in another point $y$. 

Take a general point $P=(u,v,w)\in 
\triplet {a}{b-1}{c-1}$ and consider the projection $\pi$ of $X$ from 
$L= u\otimes B\otimes C$. The image of the projection
is $Y=\PP(A)\otimes\PP(B)\otimes\PP(C)$
where $A\subset A'$ has codimension $1$. 
Furthermore, by the generality of 
$P$, $L$ does not meet $\Lambda$, as well as 
any line spanned by $y,x_i$.
It follows that the span of the tangent spaces to $Y$
at the general points $\pi(x_1), \dots , \pi(x_k)$
and containing the spaces 
$A\otimes\pi(u_i)$, $B\otimes\pi(v_i)$, 
$C\otimes\pi(w_i)$ is
also tangent in another point $\pi(y)$. Thus $Y$ is $(k,p,q,r)$-weakly defective.
By induction, we get a contradiction. 
\end{proof}

Now we are ready to state and prove our inductive criterion.

Let  $X'=\PP(A')\times\PP(B)\times\PP(C)$, 
$X''=\PP(A'')\times\PP(B)\times\PP(C)$.
Note that $A\otimes B\otimes C=\left(A'\otimes B\otimes C\right)\oplus 
\left(A''\otimes B\otimes C\right)$.
Denote by $\pi'$ and $\pi''$ the two projections.

\begin{Prop}\label{induz} ({\bf Inductive Step.})
Assume that $X'$ is $(k_1,p+k_2,q_1,r_1)$-not weakly defective 
and $X''$ is $(k_2,p+k_1,q_2,r_2)$-not weakly defective.
Then $X$ is $(k_1+k_2,p,q_1+q_2,r_1+r_2)$-not weakly defective.
\end{Prop}

\begin{proof} We specialize $k_1+k_2$ points on $X$ in order
that $k_1$ of them belong to $X_1$ and $k_2$ of them belong to $X_2$ .
Let $x_1,\ldots x_{k_1}\in X'$ and $y_1,\ldots y_{k_2}\in X''$.

Let $A\otimes\tilde v_i\otimes\tilde w_i$ for 
$i=1,\ldots p$, be subspaces.

We specialize $q_1+q_2$ points in $\PP(A)\times\PP(C)$ 
in order that  the first $q_1$ of them belong to $\PP(A')\times\PP(C)$ 
and the last $q_2$ of them belong to $\PP(A'')\times\PP(C)$ .
Call $Q_1$ the span of the first $q_1$ spaces $Bv_i$ and 
$Q_2$ the span  of the last $q_2$ spaces $Bv_i$.

We specialize $r_1+r_2$ points in $\PP(A)\times\PP(B)$ in order that 
the first $r_1$ of them belong to $\PP(A')\times\PP(B)$ and the last 
$r_2$ of them belong to $\PP(A'')\times\PP(B)$.
Call $R_1$ the span of the first $r_1$ spaces $Cw_i$ 
and $R_2$ the span  of the last $r_2$ spaces $Cw_i$.

We want to prove that $T=T_{x_1}X+\ldots +T_{x_{k_1}}X+T_{y_1}X+
\ldots +T_{y_{k_2}}X+Q_1+Q_2+R_1+R_2+A\otimes\tilde v_1\otimes
\tilde w_1+\dots +A\otimes\tilde v_p\otimes
\tilde w_p$, 
is tangent to $X$ only at $x_1,\ldots x_{k_1},y_1,\ldots y_{k_2}$.

Let $T_xX\subset T$, with $x=u\otimes v\otimes w$.
Then $\pi_1(T_xX)\subset\pi_1(T)$.
Let $u=u'+u''$, $y_j=u''_j\otimes v''_j\otimes w''_j$, $j=1,\dots,k_2$.
At least one among $u'$ and $u''$ is non zero, so let's 
assume $u'\neq 0$. Then we get
$\pi_1(T_xX)=A'\otimes v\otimes w+u'\otimes B\otimes w+u'\otimes 
v\otimes C$ while $\pi_1(T)=T_{x_1}X'+\ldots +T_{x_{k_1}}X' 
+A'\otimes v''_j\otimes w''_j+\ldots + A'\otimes\tilde v_i
\otimes \tilde w_i+ \ldots+Q_1+R_1$ (with $i=1,\dots,p$).
By the assumption that $X'$ is $(k_1,p+k_2,q_1,r_1)$-not weakly 
defective it follows that $u'\otimes v\otimes w$ is one among $x_i$.

If also $u'' \neq 0$ the same argument shows that $u''\otimes v\otimes w$ 
is one among $y_i$, which is a contradiction. Then $u''=0$, 
that is $x=u'\otimes v\otimes w$ is one among $x_i$.
 It follows that $X$ is $(k_1+k_2,p,q_1+q_2,r_1+r_2)$-not 
 weakly defective, as we wanted. 
\end{proof} 

The inductive procedure stops eventually
when we find some condition on weak  defectivity, which
does not hold. This does not means, in general, that 
our starting example was not $k$-identifiable, but merely
that we specialized the points too much, in order
to expect a meaningful answer.
\vskip 0.8cm

{\it Proof of Theorem \ref{main2} in cases $a=8, 9, 10$.}

In case $a=8$ we start with $22$ points and we want to apply iteratively the
Proposition \ref{induz}.
Splitting one $8$-dimensional vector space of the product
in a direct sum of two $4$-dimensional spaces, one sees that the
$(22,0,0,0)$-not weak defectivity of $\triplet777$
follows if one knows that $\triplet377$ is 
$(11,11,0,0)$-not weakly defective. Repeating the procedure with
the second factor, everything  reduces to prove
that $\triplet337$ is $(5,7,6,0)$-not weakly defective and
$(6,4,5,0)$-not weakly defective.
The first statement reduces to show that 
$\triplet333$ is $(3,3,3,2)$-not weakly defective and
$(2,4,3,3)$-not weakly defective. These statements have finally a reasonable size
and can be checked with a random choice of points with our Macaulay2 algorithm.
The last statement  reduces to show that 
$\triplet333$ is $(3,2,3,3)$-not weakly defective and
$(3,2,2,3)$-not weakly defective, which follows from the above check and by the Remark
\ref{ovv} (a).

In the case $a=9$ we start with $27$ points and we split the $9$ dimensional space in {\it three} $3$-dimensional summands. The inductive step is better explained by the following table 
$$
\begin{array}{ccccccccc}
a & b & c   & & k & & p &  q & r \\
9&9&9&&27&&0&0&0\\
3&9&9&&9&&18&0&0\\
3&3&9&&3&&6&6&0\\
3&3&3&&1&&2&2&2\\
\end{array}$$
The last statement can be checked again with Macaulay2.

The $a=10$ case starts as follows
$$
\begin{array}{ccccccccc}
a & b & c   & & k & & p &  q & r \\
10&10&10&&32&&0&0&0\\
5&10&10&&16&&16&0&0\\
5&5&10&&8&&8&8&0\\
5&5&5&&4&&4&4&4\\
\end{array}$$

The second statement  reduces to show that 
$\triplet144$ is $(1,7,2,2)$-not weakly defective 
and
$\triplet244$ is $(3,5,2,2)$-not weakly defective. 
Both these statements can be checked  with Macaulay2.
This concludes the proof.

\section{Proof of Theorem \ref{main}}\label{pfs}

In order to use the inductive step, we need a starting point.

\begin{Lem}\label{start}
$\triplet 111$ is $(1,0,0,0)$-not and $(0,1,1,1)$-not weakly defective.
\end{Lem}
\begin{proof}
The first fact is true for any smooth variety,  see  Example \ref{nonbitangent}.
For the second one, we consider $X=\PP(A)\times\PP(B)\times\PP(C)$ where $A$, $B$, $C$ have all dimension $2$
and we choose basis $A=\langle a_0, a_1\rangle$, $B=\langle b_0, b_1\rangle$, $C=\langle c_0, c_1\rangle$.
Then, without loss of generality, we may consider the span $T=A\otimes b_0\otimes c_0+a_0\otimes B\otimes c_1+a_1\otimes b_1\otimes C$.
In the monomial basis of $A\otimes B\otimes C$ this span contains all the monomials with the only exception of $a_0\otimes b_1\otimes c_0$
and $a_1\otimes b_0\otimes c_1$. Then, a vector $v=\sum x_{ijk}a_i\otimes b_j\otimes c_k$,
belongs to $X\cap \PP(T)$ if all the $2\times 2$-minors of the two following flattening matrices vanish
$$\left[\begin{array}{rrrr}
x_{000}&x_{001}&x_{100}&0\\
0&x_{011}&x_{110}&x_{111}
\end{array}\right]\qquad\qquad
\left[\begin{array}{rrrr}
x_{000}&0&x_{100}&x_{110}\\
x_{001}&x_{011}&0&x_{111}
\end{array}\right]$$
A straightforward check on the minors shows that $X\cap \PP(T)$ consists of the following 
six lines in the $5$-dimensional space $\PP(T)=\{x_{010}=x_{101}=0\}$
(Macaulay2 can be helpful at this step) 

$r_0=V(x_{001},x_{000},x_{100},x_{110})$

$r_1=V(x_{000},x_{100},x_{110},x_{111})$

$r_2=V(x_{100},x_{110},x_{111},x_{011})$

$r_3=V(x_{110},x_{111},x_{011},x_{001})$

$r_4=V(x_{111},x_{011},x_{001},x_{000})$

$r_5=V(x_{011},x_{001},x_{000},x_{100})$

which have the property that, for $i\neq j$ 
$$r_i\cap r_j=\left\{\begin{array}{cl}
\textrm{one point} &\textrm{if\ }i=j+1, j-1\textrm{\ mod\ }6\\
\emptyset&\textrm{otherwise}
\end{array}\right.$$
It follows that $\PP(T)$ is not tangent anywhere,
because the tangent space at a point meets $X$ in three concurrent lines. 
This proves that $X$ is  $(0,1,1,1)$-not weakly defective.
\end{proof}

\begin{Rmk}
We will use affine spaces whose dimension is a power of $2$, as
well as  sets of points or subspaces whose number is 
expressed in terms of powers of $2$, essentially because they allow
the following recursive application of Lemma \ref{induz}:

Assume we want to prove that
$\triplet{2^\alpha-1}{2^\beta-1}{2^\gamma-1}=\PP(A)\otimes
\PP(B)\otimes\PP(C)$ is $(2x,2^u, 2^v, 2^w)$-not 
weakly defective.  Then, by splitting the first linear space $A$ in a 
direct sum of two subspaces of dimension $2^{\alpha-1}$ and balancing 
the splitting of the number of points and linear spaces, by 
Proposition \ref{induz}  it is sufficient to prove that
$\triplet{2^{\alpha-1}-1}{2^\beta-1}{2^\gamma-1}$ is 
$(x,2^u, 2^{v-1}, 2^{w-1})$-not weakly defective. 

We will use this trick so often, in the arguments below.
\end{Rmk}

The final statement will be that, if we order the dimensions
so that $1\leq \alpha\leq \beta\leq \gamma$, then 
$X=\triplet{2^\alpha-1}{2^\beta-1}{2^\gamma-1}$ is $(k,0,0,0)$-not weakly defective,
for $k\leq  2^{\alpha+\beta-2}.$

Before showing this fact, we need a series of lemmas.

\begin{Prop} Assume that $X=\triplet{2^\alpha-1}{2^\beta-1}{2^\gamma-1}$ is
not $(k,0,0,0)$-not weakly defective. 
Then also $X'=\triplet{2^\alpha-1}{2^\beta-1}{2^{\gamma}}$ is $(k,0,0,0)$-not 
weakly defective.
\end{Prop}
\begin{proof} 
By Lemma \ref{riduz}.
\end{proof}

So, in order to prove  Theorem \ref{main}, we can reduce ourselves to 
the case $\beta=\gamma$, $k=2^{\alpha+\beta-2}$.

\begin{Lem}\label{zero} Take $X=\triplet{2^{a_1}-1}{2^{a_2}-1}{2^{a_3}-1}$, with 
$a_1,a_2,a_3\geq 1$. Pick non-negative integers $u_1,u_2, u_3$
such that $u_i\leq a_j+a_k-2$, whenever $\{i,j,k\}=\{1,2,3\}$.
Then $X$ is $(0,2^{u_1},2^{u_2},2^{u_3})$-not weakly defective. 
\end{Lem}
\begin{proof} We make induction on the sum $a_1+a_2+a_3$.

If $a_1=a_2=a_3=1$, then the numerical conditions imply that
$u_1=u_2=u_3=0$ and the conclusion follows from the fact that
$\triplet 111$ is $(0,1,1,1)$-not weakly defective, which holds by 
 Lemma \ref{start}. 

Assume $a_1>1$ and split the first projective space in a sum of 
two spaces of dimension $2^{a_1-1}$. Then there are three
possibilities:

(1) Assume $u_2=u_3=0$. Then, by using Lemma \ref{induz}, the claim
reduces to prove that $\triplet{2^{a_1-1}-1}{2^{a_2}-1}{2^{a_3}-1}$ is 
 $(0,2^{u_1}, 1, 1)$-not weakly defective and it is 
 $(0,2^{u_1}, 0, 0)$-not weakly defective. The second condition 
 is contained in the first. Since $a_1>1$, the six numbers 
 $a_1-1, a_2, a_3, u_1, 0 , 0$ fulfill the numerical inequalities
of the statement. Hence the claim follows by induction, in this case.

(2) Assume $u_3>u_2=0$. Then the claim
reduces to prove that $\triplet{2^{a_1-1}-1}{2^{a_2}-1}{2^{a_3}-1}$ is 
 $(0,2^{u_1}, 1, 2^{u_3-1})$-not weakly defective and it is 
 $(0,2^{u_1}, 0, 2^{u_3-1})$-not weakly defective. The second condition 
 is contained in the first. One checks that the six numbers 
 $a_1-1, a_2, a_3, u_1, 0 , u_3-1$ fulfill the numerical inequalities
of the statement. Hence the claim follows by induction.

(3) Assume $u_2,u_3>0$. Then the claim
reduces to prove that $\triplet{2^{a_1-1}-1}{2^{a_2}-1}{2^{a_3}-1}$ is
$(0,2^{u_1}, 2^{u_2-1}, 2^{u_3-1})$-not weakly defective. 
One checks that the six numbers 
 $a_1-1, a_2, a_3, u_1, u_2-1 , u_3-1$ fulfill the numerical inequalities
of the statement. Hence the claim follows by induction.
\end{proof}

\begin{Lem}\label{uno} Take $X=\triplet{2^{a_1}-1}{2^{a_2}-1}{2^{a_3}-1}$, with 
$a_1,a_2,a_3\geq 1$. Pick non-negative integers $u_1,u_2, u_3$
such that $u_i\leq a_j+a_k-2$, whenever $\{i,j,k\}=\{1,2,3\}$.
Then $X$ is $(1,2^{u_1}-1,2^{u_2}-1,2^{u_3}-1)$-not weakly defective. 
\end{Lem}
\begin{proof} We make induction on the sum $a_1+a_2+a_3$.

If $a_1=a_2=a_3=1$, then the numerical conditions imply that
$u_1=u_2=u_3=0$ and the conclusion follows from the fact that
$\triplet 111$ is $(1,0,0,0)$-not weakly defective (Lemma \ref{start}).

Assume $a_1>1$ and split the first projective space in a sum of 
two spaces of dimension $2^{a_1-1}$. Then there are three
possibilities:

(1) Assume $u_2=u_3=0$. Then, by using Lemma \ref{induz}, the claim
reduces to prove that $\triplet{2^{a_1-1}-1}{2^{a_2}-1}{2^{a_3}-1}$ is 
 $(0,2^{u_1}, 1, 1)$-not weakly defective and it is 
 $(1,2^{u_1}-1, 0, 0)$-not weakly defective. The first condition follows 
 by the previous Lemma \ref{zero}. For the second condition, 
notice that since $a_1>1$, the the six numbers 
 $a_1-1, a_2, a_3, u_1, 0 , 0$ fulfill the numerical inequalities
of the statement (and $0=2^0-1$). 
Hence the claim follows by induction, in this case.

(2) Assume $u_3>u_2=0$. Then the claim
reduces to prove that $\triplet{2^{a_1-1}-1}{2^{a_2}-1}{2^{a_3}-1}$ is
 $(0,2^{u_1}, 1, 2^{u_3-1})$-not weakly defective and it is
 $(1,2^{u_1}-1, 0, 2^{u_3-1}-1)$-not weakly defective. 
The first condition follows by the previous Lemma. The second condition 
follows by induction, since one checks that the six numbers 
 $a_1-1, a_2, a_3, u_1, 0 , u_3-1$ fulfill the numerical inequalities
of the statement.

(3) $u_2,u_3>0$. Then the claim
reduces to prove that $\triplet{2^{a_1-1}-1}{2^{a_2}-1}{2^{a_3}-1}$ is
$(0,2^{u_1}, 2^{u_2-1}, 2^{u_3-1})$-not weakly defective and it is 
$(1,2^{u_1}-1, 2^{u_2-1}-1, 2^{u_3-1}-1)$-not weakly defective. 
One checks that the numerical conditions in the statement are still 
fulfilled, by the six numbers $a_1-1, a_2, a_3, u_1, u_2-1 , u_3-1$. 
Hence the claim follows by induction.
\end{proof}

Now we are ready to prove:

\begin{Thm}\label{premain}
 $X=\triplet {2^\alpha-1}{2^\beta-1}{2^\beta-1}$ is $(k,0,0,0)$-not 
 weakly defective,  for $k\leq 2^{\alpha+\beta-2}$.
\end{Thm}
\begin{proof} 
Write $\alpha+\beta-2=2p+e$, where $e$ is the remainder. 

Now we start our reduction.

($A_1$) One can split the
vector space in the middle as a sum of two spaces of
dimension $2^{\beta-1}$. By using Lemma \ref{induz}, it turns out that
$X$ is $(2^{\alpha+\beta-2},0,0,0)$-not weakly defective
when $\triplet {2^\alpha-1}{2^{\beta-1}-1}{2^\beta-1}$ is 
 $(2^{\alpha+\beta-3},0,2^{\alpha+\beta-3},0)$-not weakly defective.

($A_2$) Splitting now the third vector space as a sum of
two spaces of dimension $2^{\beta-1}$, and using the Lemma, 
this reduces to prove that

$\triplet {2^\alpha-1}{2^{\beta-1}-1}{2^{\beta-1}-1}$ is 
$(2^{\alpha+\beta-4},0,2^{\alpha+\beta-4},2^{\alpha+\beta-4})$-not 
weakly defective.

($A_3$) Now repeat the procedure, splitting the space in the middle:
Everything  reduces to prove that 

$\triplet {2^\alpha-1}{2^{\beta-2}-1}{2^{\beta-1}-1}$ is 
$(2^{\alpha+\beta-5},0,2^{\alpha+\beta-4}+2^{\alpha+\beta-5},
2^{\alpha+\beta-5})$-not weakly defective.

Now split again the third vector space, and
repeat the steps. At the end of the $(\alpha+\beta-2)$-th step, after
the computation, we find out that we need just to prove that

 $\triplet{2^\alpha-1}{2^{\beta-p-e-1}}{2^{\beta-p-1}}$ is $(1,0,
 \sum_{i=0}^{p+e-1} 2^i, \sum_{i=0}^{p-1} 2^i)$-not weakly defective.

Notice that all these steps can be performed, because $\beta-p\geq 
\beta-p-e\geq 1$. Indeed  we have $\alpha\leq \beta$, thus 
$2\beta-2\geq 2p+e$, hence  $2\beta\geq 2p+e+2>2p+2e$.

Now, $\sum_{i=0}^{p+e-1} 2^i=2^{p+e}-1$ while $\sum_{i=0}^{p-1} 2^i=2^p-1$.
Moreover
\begin{gather} \nonumber
p+e \leq \alpha+(\beta-p)-2 \qquad \mbox{ since } 2p+e=\alpha+\beta-2 \\ 
\nonumber p\leq \alpha+(\beta-p-e)-2 \qquad \mbox{ since } 
2p = \alpha+\beta-e-2. \end{gather} 

Thus we may apply Lemma \ref{uno}, and see that 
$\triplet{2^\alpha-1}{2^{\beta-p-e-1}}{2^{\beta-p-1}}$ is $(1,0,2^{p+e}-1, 2^p-1)$-not 
weakly  defective. The result is settled.
\end{proof}

When $\alpha=\beta$, i.e. when the product is balanced, we find that $X$
is $k$-identifiable for $k\leq 2^{2\alpha-2}$.
\smallskip

\noindent{\it Proof of Theorem \ref{main}} Fix $\alpha,\beta$ 
maximal such that  $2^\alpha\leq a$ and $2^\beta\leq b$. 
Then, by the previous Theorem,  $\triplet{2^\alpha-1}{2^\beta-1}{2^\beta-1}$ 
is $(k,0,0,0)$-not weakly defective, for $k\leq 2^{\alpha+\beta-2}=
2^\alpha2^\beta/4$. Thus also $\PP(A)\otimes\PP(B)\otimes\PP(C)=
\triplet {a-1}{b-1}{c-1}$ is $(k,0,0,0)$-not 
weakly defective, for  $k\leq 2^\alpha2^\beta/4$.
The conclusion follows. 
\hfill\qed
\vskip 0.5cm

Comparing our result with the maximal $k$ for which the 
identifiability of $\PP(A)\otimes\PP(B)\otimes\PP(C)$ makes sense, i.e.
$$k_{max}=   \lfloor  \frac{abc}{a+b+c-2} \rfloor.$$
(see Proposition \ref{maxk}) and considering that 
$ab/3\leq k_{max}\leq ab$,
we see that the bound in the Theorem is, at least log-asymptotically,
sharp, as explained in the Introduction.

In any events, it improves Kruskal's bound for identifiability.

\begin{Rmk} In principle, there are no obstructions in repeating the
argument of Theorem \ref{premain}, when we substitute powers of $2$
with powers of any other integer $p>1$. The final statement is:

\noindent{\it $X=\triplet{p^\alpha-1}{p^\beta-1}{p^\beta-1}$ is 
$(k,0,0,0)$-not weakly defective, for $k\leq p^{\alpha+\beta-2}$.}

The proof is achieved very similarly, by splitting, step by step,
a vector space of dimension $p^n$ into $p$ spaces of dimension $p^{n-1}.$
(see e.g. the case $a=9$ in the proof of Theorem \ref{main2}).

We can use this statement, instead of
Theorem \ref{premain}, in the proof of Theorem \ref{main},
obtaining another bound which implies $k$-identifiability.
 
In most cases, however, the new bound is weaker than the one of the Theorem
\ref{main}. On the other hand, in some specific case, typically 
when powers of $3$ are involved, it can be stronger.

To give an example, let us consider $X=\triplet{26}{26}{26}$.
Using Theorem \ref{main}, we obtain $k$-identifiability for
$k\leq 2^{4+4-2}=64$. Using powers of $p=3$, instead, we
get $k$-identifiability for $k\leq 3^{3+3-2}=81$. It is an improvement,
but still a long way from $k_{max}=249$.
\end{Rmk}

\section{Some examples in low dimension}
\label{fin}

In this section, we study the $k$-identifiability of Segre products
$X=\PP(A)\otimes\PP(B)\otimes\PP(C)$, when the dimensions $a,b,c$ are small. 
We also provide a proof for Theorem \ref{main3}.

\smallskip

\noindent{\it Proof of Theorem \ref{main3}}.
Consider $X=\triplet333$. This product is
$5$-identifiable, by Kruskal's criterion.
On the other hand, accordingly with Proposition \ref{maxk},
one may ask about the $6$-identifiability of $X$.

We are able to prove that $X$ is {\it not} $6$-identifiable,
and the general point in $S^6(X)$ sits in exactly two $6$-secant,
$5$-planes. From the tensorial point of view,
this means that a general $4\times 4 \times 4$ tensor of rank $6$,
can be written as a sum of $6$ decomposable tensors in exactly
$2$ ways (up to scalar multiplication and permutations).

The reason relies in the fact that through $6$ general points 
$x_1,\dots,x_6$ of $X=\triplet{3}{3}{3}$, one can draw an elliptic normal 
curve $\Gamma$ of degree $12$, which spans a projective space $L=\PP^{11}$,
containing the linear span of $x_1,\dots,x_6$. So, a general point
$u\in S^6(X)$ lies in a linear space $L$ spanned by an elliptic normal
curve $\Gamma\subset X$. By \cite{CC2}, Proposition 5.2, it is known
that $\Gamma$ has $6$-secant order $2$, i.e. there are exactly
two $5$-planes, $6$-secant to $\Gamma$, inside $L$. By \cite{CC2}
Proposition 2.4, if we prove that $\Gamma$ coincides with the contact
locus of a general $6$-tangent hyperplane,
also $X$ must have $6$-secant order equal to $2$.
This last fact can be checked by our Macaulay2 algorithm.
Unfortunately, the existence of an elliptic normal curve, of degree $12$,
 passing through $6$ randomly chosen points of $X$, gives only
a probabilistic argument for the existence of such a curve passing through
$6$ general points of $X$. To overcome this problem, 
we offer the following theoretical  argument.

 Consider the projections $z_{i1},\dots,z_{i6}$ of $x_1,\dots,x_6$, 
into the $i$-th copy of $\PP^3$, so that $z_{i1},\dots,z_{i6}$ are general
points of $\PP^3$. Normal elliptic curves $C$ passing through the
$6$ points of $\PP^3$ are given by pairs of quadrics through
 the points, so they are parametrized by the Grassmannian $G$ of lines 
 in the space $\PP^3$ of quadrics through $z_{i1},\dots,z_{i6}$.
In order that three normal elliptic curves $C,C',C''$ 
in the three copies of $\PP^3$,
correspond to the same abstract curve, they need to 
differ by an element of $PGL(3)$. 
So, once we have $C$ ($4$ parameters), we can choose 
 $\phi,\psi\in PGL(3)$ for the two remaining maps $C\to\PP^3$
 (thus a total of $4+15+15=34$ parameters). On the other hand,
 we need to impose that $\phi(C)=C'$ (resp $\psi(C)=C''$)
 pass through $z_{21},\dots,z_{26}$ (resp. $z_{31},\dots,z_{36}$).
 Since each point imposes $2$ conditions,
we get a total of $24$ algebraic conditions on the $34$ parameters.

Moreover, if we want that after
this correspondence, $C,C',C''$ are projection of the same curve
passing through $x_1,\dots,x_6$,  we also need that the projectivity
$\phi: C\to C'$ (resp. $\psi: C\to C''$), composed with 
the automorphisms of the curves which sends $z_{11}$ to $z_{21}$
(resp. $z_{11}$ to $z_{31}$), also sends any $z_{1i}$ to $z_{2i}$
(resp. $z_{1i}$ to $z_{3i}$), for $i\geq 2$. 
This gives $10$ more conditions,
which are algebraic on the coefficients 
 of the two quadrics and the entries of the matrices of $\phi,\psi$.
 
So, we have a total of $34$ conditions, which are algebraic
on the $34$ parameters, i.e. on the projective coordinates
of $G\times PGL(3)\times PGL(3)$. Thus we get at least
a finite number of curves passing through $x_1,\dots,x_6$, 
for a general choice of the points.  
\hfill\qed
\vskip .5cm

\begin{Rmk} In the previous example, notice that when 
the three projections of the points $x_1,\dots, x_k$ differ
by a projectivity, then the number of conditions decreases,
and we find infinitely many normal elliptic curves.

It is easy to see that this implies that a point in the
secant variety $S_6$ of any of these curves, belongs indeed 
to infinitely many $6$-secant spaces.
\end{Rmk}

The case of products of projective spaces of dimension $3$ is
 particularly interesting, due to its applications to 
statistical studies on DNA strings.

If we have many substrings of DNA strings, each formed by three positions,
and we record the occurrence of the four bases in each position, we get a
distribution which can be  arranged in a $4\times 4\times 4$ tensor $T$.
The rank $k$ of $T$ suggests the existence of $k$ different types
of substrings, in the probe, such that for each type, the distribution of
bases is independent. So $T$ is the sum of $k$ tensors
$T_1,\dots,T_k$, of rank $1$.

An obvious question concerns the possibility of recovering the
$k$ tensors $T_i$, starting from $T$. When $k\geq 7$, this possibility
is excluded, since $7$ exceeds the maximum given in Proposition \ref{maxk}.
For $k\leq 5$, $k$-identifiability (by Kruskal's criterion)
tells us that, at least theoretically, the reconstruction is possible.

The amazing situation happens for $k=6$. Although one could
expect $6$-identifiability, Theorem \ref{main3} shows that there
are exactly two sets of tensors of rank $1$, whose sum is $T$. 
Hence, at least over the complex field, there are exactly  two
different sets of distributions, in the $6$ types, that produce
the same distribution $T$. 
 \medskip

In \cite{AOP}  6.3 one finds the list of known 
Segre varieties $X=\PP(A)\times\PP(B)\times\PP(C)=
\triplet {a-1}{b-1}{c-1}$ (with $a\le b\le c$) such that
the dimension of $k$-th secant variety  is smaller than 
the expected value. 
Recall that when the dimension of $S_k(X)$ is smaller than the expected value,
i.e. when the variety $X$ is $k$-defective,
then the $k$-identifiability necessarily fails.

A list of known Segre varieties $X$ 
which are not $k$-identifiable, i.e. such that
the general tensor of rank $k$ in $\CC^a\otimes\CC^b\otimes\CC^c$
has not a unique decomposition, is the following
(for $k<k_{max}$):

{\small
$$\begin{array}{lccl}
&(a,b,c)&k& notes\\
\hline\\
\textrm{defective}&c\ge (a-1)(b-1)+3&(a-1)(b-1)+2\le k 
&\textrm{\cite{AOP}}\\
{\textrm{unbalanced}}&&k<\min \left(c,ab\right)&\\
\hline\\
\textrm{defective}&(3,4,4)&5&\textrm{\cite{AOP}}\\
\hline\\
\textrm{defective}&(3,b,b)\quad b\textrm{\ odd}&
\frac{3b-1}{2}&\textrm{\cite{S}}\\[0.1cm]
\hline\\
{\textrm{w. defective}}&3\le a &(a-1)(b-1)+1& \binom d{(a-1)(b-1)+1}
\\ \textrm{unbalanced}& c\ge (a-1)(b-1)+2 &&\textrm{\ decompositions} \\
\textrm{}&&&\textrm{where\ } d={{a+b-2}\choose {a-1}}\\
\textrm{}&&& \textrm{(Theorem \ref{unbal})}\\
\hline\\
{\textrm{w. defective}}&(4,4,4)&6&2\textrm{\ decompositions} \\ \textrm{}&&& \textrm{(Theorem \ref{main3})}\\
\hline\\ 
{\textrm{w. defective}}&(3,6,6)&8& (**)
\end{array}$$
}

A computer check shows that this list is complete for $c\le 7$.
 In the last case marked with (**), the contact variety is a
$4$-fold in $\PP^{39}$ of degree $108$.
This case needs an ''ad hoc'' analysis which goes beyond the space
of the present note and will be addressed in a forthcoming paper \cite{CMO}.

\vskip.5cm

In the unbalanced case, the identifiability can be proved theoretically.

\begin{Prop}\label{rankab} The general tensor of rank $(a-1)(b-1)$ in
$\PP(\CC^a\otimes\CC^b\otimes\CC^c)$
has a unique decomposition as sum of $(a-1)(b-1)$ summands in
$\triplet {a-1}{b-1}{c-1}$ for $c\ge (a-1)(b-1)$.
\end{Prop}
\begin{proof}
Let $\phi\in\CC^a\otimes\CC^b\otimes\CC^c$ be general of rank $(a-1)(b-1)$.
It induces the flattening  contraction operator 
$$A_{\phi}\colon(\CC^c)^{\vee}\to\CC^a\otimes\CC^b$$
which has still rank $(a-1)(b-1)$, by the assumption $c\ge (a-1)(b-1)$. 
Indeed, if $\phi=\sum_{i=1}^{(a-1)(b-1)}u_i\otimes v_i\otimes w_i$
with $u_i\in \CC^a, v_i\in\CC^b, w_i\in\CC^c$, 
where $w_i$ can be chosen as part of a basis of $C$,
 then $\textrm{Im\ }A_{\phi}$ is the span of the representatives
of $v_i\otimes w_i$ for $i=1,\ldots , (a-1)(b-1)$.
It is well known that the projectification of this span,
whose dimension is smaller than the codimension of the Segre variety
$Y=\PP^{a-1}\times\PP^{b-1}\subset\PP(\CC^a\otimes\CC^b)$, meets $Y$
only in these $(a-1)(b-1)$ points
(see for example the Theorem 2.6 in \cite{CC1}).
The claim follows.
\end{proof}

\begin{Prop}
When $c=(a-1)(b-1)$ or $c=(a-1)(b-1)+1$, then the rank of a generic tensor in 
 $\PP(\CC^a\otimes\CC^b\otimes\CC^c)$  is $ab-a-b+2$.
\end{Prop}
\begin{proof}
When $c\ge (a-1)(b-1)+1$, we are in the unbalanced case, 
according to the definition 4.2 of \cite{AOP}. In this case the 
generic rank is $\min\{c, ab\}$ by (ii) 
of the Theorem 4.4  of \cite{AOP}. 

Assume $c=(a-1)(b-1)$. Using the same technique, we show that the 
secant variety $S_k(\triplet {a-1}{b-1}{c-1})$ has the expected dimension, 
for $k\le (a-1)(b-1)$, and fills the  ambient space, for $k=(ab-a-b+2)$. 

Indeed, with the notations of \cite{AOP},
$T(a-1,b-1,ab-a-b;(a-1)(b-1);0,0,0)$ reduces to $T(a-1,b-1,0;1;0,0,ab-a-b)$, 
which is true and subabundant, while $T(a-1,b-1,ab-a-b;ab-a-b+2;0,0,0)$ reduces 
(for $b\ge 3$) to $T(a-1,b-1,0;1;0,0,ab-a-b+1)$ and 
$T(a-1,b-1,0;2;0,0,ab-a-b-1)$ which are both superabundant and true.
\end{proof}

\begin{Prop}\label{different}
Assume $c\ge (a-1)(b-1)+2$. Then  the generic rank in  
$\PP(\CC^a\otimes\CC^b\otimes\CC^c)$ is at least $(a-1)(b-1)+2$,
and it is equal to $(a-1)(b-1)+2$ in the border case $c=(a-1)(b-1)+2$.
The number of different decomposition of a general 
tensor of rank $(a-1)(b-1)+1$ is ${d\choose {(a-1)(b-1)+1}}$
where $d=\deg (\PP^{a-1}\times\PP^{b-1})={{a+b-2}\choose {a-1}}$. 
This number is always bigger than $1$, with the only exception
$a=b=2$.
\end{Prop}
\begin{proof}
We apply the same argument of the proof of the Proposition \ref{rankab}.
The unique difference is that, now, the dimension of the
projectification of $\textrm{Im\ }A_{\phi}$ equals the codimension
of $\PP^{a-1}\times\PP^{b-1}$. Thus we get $d$ points of intersection.
Any choice of $(a-1)(b-1)+1$ among these $d$ points, 
yields a decomposition.
\end{proof}

\begin{Rmk}
The case $a=b=3$ of Prop. \ref{different} is connected to the work of tenBerge,
who showed in \cite{tB},
that there are six different decompositions of a general rank $5$
tensor in $\CC^3\otimes\CC^3\otimes \CC^5$, chosen taking $5$ among $6$ possible summands. Our argument,
which we gave for $c\ge 6$, can be extended to the case
$c=5$ and $k=k_{max}=5$, and it gives a geometric explanation of this phenomenon,
indeed the six possible summands correspond to the six intersection points of $\PP^2\times\PP^2$ 
with a general $\PP^4$.
\end{Rmk}

As a consequence of the two previous results, we get:

\begin{Thm}\label{unbal} 
Assume $c\ge (a-1)(b-1)+2$. Then the general tensor of rank $k$ in
$\PP(\CC^a\otimes\CC^b\otimes\CC^c)$
has a unique decomposition as sum of $k$ summands in
$\triplet {a-1}{b-1}{c-1}$ if and only if $k\le (a-1)(b-1)$.
\end{Thm}

\section{Products with many factors}\label{many}

At the cost of the growth of the notation,
we can generalize the statement of our main Theorem
\ref{main}, to products of many vector spaces.

In this section, we simply list the corresponding
definitions and results. The proofs are absolutely
straightforward, following the pattern
of the corresponding arguments in the previous sections.
Only the initial step of the induction needs an extra
argument, which is displayed in  Lemma \ref{zerostep} below.
\smallskip

For a given set of complex vector spaces $A_1,\dots,A_n$,
with $n\geq 3$ and $\dim A_i\geq 2$, let us give the general:

\begin{Def}\label{p_i}  A Segre product
 $X=\PP(A_1)\times\dots\times\PP(A_n)$ 
is called {\it $(k,p_1,\dots,p_n)$-not weakly defective} if:
 
for $k$ general points $x_1,\ldots x_k\in X$, 

for $p_i$ general points $w_{ij}\in\PP(A_1)\times\dots\times\hat\PP(A_i)
\times\dots\times\PP(A_n)$,

the span of the spaces $T_{x_i}X$, $A_i\otimes w_{ij}$ 
contains  $T_xX$ if and only if $x=x_i$, for some $i=1,\ldots k$.
Otherwise $X$ is called {\it $(k,p_1,\dots,p_n)$-weakly defective}.
\end{Def}

\begin{Rmk}\label{ovvn} (a) With the previous notation, by semicontinuity
it is clear that when $X$ is  $(k,p_1,\dots,p_n)$-not weakly defective,
then it is also  $(k',p'_1,\dots,p'_n)$-not weakly defective,
whenever $(k',p'_1,\dots,p'_n)\leq (k,p_1,\dots,p_n)$, in the strict ordering.

(b) By semicontinuity, $X$ is  $(k,p_1,\dots,p_n)$-not weakly defective
whenever one gets that for {\it particular} sets of points
$\{x_i\}$, $\{w_{ij}\}$, as above,
then the span of 
$T_{x_i}X$ and all $A_i\otimes w_{ij}$
contains  $T_xX$ if only if $x=x_i$, for some $i=1,\ldots k$.
 
(c) By Proposition \ref{criterio}, one gets soon that
 $(k,0,\dots,0)$-not weakly defective implies  $k$-identifiable.
\end{Rmk}

\begin{Lem} \label{riduzn} 
Consider $X=\PP(A_1)\times\dots\times\PP(A_n)$ and assume that,
for a choice of subspaces $A'_i\subset A_i$,
the product  $\PP (A'_1)\times\dots\times\PP(A'_n)$ is
$(k,p_1,\dots,p_n)$-not weakly defective. Then $X$ is 
$(k,p_1,\dots,p_n)$-not weakly defective.
\end{Lem}

The inductive criterion can be rephrased as follows,
 always following the lines in \cite{AOP}.

\begin{Prop}\label{induzn} {\bf Inductive Step}
Split the vector space $A_i$ in the sum of two
spaces $A'_i$ and $A''_i$.
Let $X'=\PP(A_1)\times\dots\times \PP(A'_i)\times\dots\times\PP(A_n)$, 
$X''=\PP(A_1)\times\dots\times \PP(A''_i)\times\dots\times\PP(A_n))$, 

Assume that the product $X'$ is $(k_1,p'_1,\dots,p_i+k_2,\dots,p'_n)$-not 
weakly defective and the product $X''$ is $(k_2,p'_1,\dots,p_i+k_1,\dots,p'_n)$-not 
weakly defective.
Then, setting $p_j=p'_j+p''_j$ for $j\neq i$, we get that $X$ is 
$(k_1+k_2,p_1,\dots,p_i, \dots, p_n)$-not weakly defective.
\end{Prop}

Now we use again the previous criterion, when the dimension of the
vector spaces are powers of $2$, i.e. when $\dim(A_i)=2^{\alpha_i}$,
for all $i$. We agree to order the spaces, so that 
$$ \alpha_1\leq \dots\leq \alpha_n. $$

The following numerical criterion is the exact generalization of Lemmas
\ref{zero} and \ref{uno}.

\begin{Lem}\label{zerostep} Take $X=\PP(A_1)\times\dots\times\PP(A_n)$, with 
$n\geq 3$ and $\dim(A_i)=2^{\alpha_i}\geq 2$. 
Pick non-negative integers $u_1,\dots, u_n$
such that, for all $i$: 
$$u_i\leq \alpha_1+\dots+\hat \alpha_i+\dots+\alpha_n-(n-1).$$
Then $X$ is  $(0,2^{u_1},\dots,2^{u_n})$-not weakly defective and
$(1,2^{u_1}-1,2^{u_2}-1,2^{u_3}-1)$-not weakly defective. 
\end{Lem}
\begin{proof} The proof goes by induction. For the inductive step,
one can follow the proof of Lemmas \ref{zero} and \ref{uno}, 
rephrased for products of many vector spaces.
Thus we only need to check the starting points of the induction, namely
that $Y_n=\PP^1\times\dots\times\PP^1$ is $(1,0,\dots,0)$-not weakly defective
and $(0,1,\dots,1)$-not weakly defective.

The first fact follows soon, as $\PP^1\times\dots\times\PP^1$ is smooth,
so that the general tangent hyperplane is not bitangent. 

The second fact follows by induction on $n$. Namely it is true for $n=3$,
as observed in Lemma \ref{uno}. For general $n$,
write $Y_n=\PP(A_1)\times\dots\times \PP(A_n)$, with
$\dim(A_i)=2$, and split $A_1$ in a direct sum
of two $1$-dimensional spaces $A'$, $A''$. 
Using Lemma \ref{riduzn}, one has thus to prove
that $Y_{n-1}=\PP^0\times\PP^1\times\dots\times\PP^1$ 
is $(0,1,0,\dots,0)$-not weakly defective
and $(0,0,1,\dots,1)$-not weakly defective.
The former claim is obvious. The latter 
follows by induction.
\end{proof}

We get:

\begin{Prop}\label{pren}
Take $X=\PP(A_1)\times\dots\times\PP(A_n)$, with 
$n\geq 3$ and $\dim(A_i)=2^{\alpha_i}\geq 2$. 
Order the $\alpha_i$'s so that $\alpha_1\leq\dots\leq \alpha_n$.
Then  $X$ is not $k$-weakly defective, for 
$k\leq 2^{\alpha_1+\dots +\alpha_{n-1}-(n-1)}$.
\end{Prop}

It follows that:

\begin{Thm} Take $X=\PP(A_1)\times\dots\times\PP(A_n)$, with 
$n\geq 3$ and $\dim(A_i)=a_i\geq 2$ and, for all $i$,
take $\alpha_i$ maximal, such that $a_i\geq 2^{\alpha_i}$. 
Then $X$ is $k$-identifiable, for 
$$k\leq 2^{\alpha_1+\dots +\alpha_{n-1}-(n-1)}.$$
\end{Thm}

Comparing our result with the maximal $k$ for which the 
identifiability of $\PP(A_1)\times\dots\times\PP(A_n)$ makes sense, 
which, in the case of a product of many factors, reads as:
$$k_{max}=\lfloor \frac{\prod_{i=1}^{n-1}a_i}
{1+\frac{\sum_{i=1}^{n-1}a_i- (n-1)}{a_n}}\rfloor 
$$
we see again that the bound in the Theorem is log-asymptotically sharp.

The inequality of the theorem can be written as
$$k\le 2^{\left(\sum_{i=1}^{n-1}\lfloor\log_2{a_i}-1\rfloor\right)}$$
Since $2^{\alpha_i}\ge\frac{a_i+1}{2}$ we get the general tensor of rank 
$k$ is $k$-identifiable if
$$k\le\frac{\prod_{i=1}^{n-1}(a_i+1)}{2^{2n-2}}$$

In \cite{SB} Kruskal bound was extended to the case of $n$ factors.
A sufficient condition for the $k$-identifiability of the general tensor of rank $k$ is  
$$2k+n-1\le\sum_{i=1}^n\min(k,a_i)$$

To compare with our condition, in the hypercubic case where $a_i=a$, the bound in \cite{SB}
is $$k\le\frac{n(a-1)+1}{2}$$
while our bound is
$$k\le 2^{(n-1)(\lfloor\log_2{a}-1\rfloor)}$$
For $a\ge 4$ we get also the weaker, but more handy, inequality
$$k\le\left(\frac{a+1}{4}\right)^{n-1}$$

\begin{Ex}
Instead of giving the proofs, which, we repeat, are analogue
to the proofs of the statement of section \ref{pfs}, let us see how
the reduction works in a concrete example.

Take $A_1=\dots =A_5=\CC^{16}$ and consider $X=\PP(A_1)\times\dots\times
\PP(A_5)$. We want to prove that $X$ is $k$-not weakly defective
for $k=2^{4+4+4+4-4}=4096$.

The reduction step starts as in the following table:

$$
\begin{array}{ccccccccccccc}
A_1 & A_2 & A_3 & A_4 & A_5  & & k & & p_1 &  p_2 & p_3 & p_4 & p_5 \\
16 & 16 & 16 & 16 & 16  & & 4096  & &    0 &    0 &   0 &   0 &   0 \\
 8 & 16 & 16 & 16 & 16  & & 2048  & & 2048 &    0 &   0 &   0 &   0 \\
 8 &  8 & 16 & 16 & 16  & & 1024  & & 1024 & 1024 &   0 &   0 &   0 \\
 8 &  8 &  8 & 16 & 16  & &  512  & &  512 &  512 & 512 &   0 &   0 \\
 8 &  8 &  8 &  8 & 16  & &  256  & &  256 &  256 & 256 & 256 &   0 \\
 8 &  8 &  8 &  8 &  8  & &  128  & &  128 &  128 & 128 & 128 & 128 \\
 4 &  8 &  8 &  8 &  8  & &   64  & &  192 &   64 &  64 &  64 &  64 \\
 4 &  4 &  8 &  8 &  8  & &   32  & &   96 &   96 &  32 &  32 &  32 \\
 4 &  4 &  4 &  8 &  8  & &   16  & &   48 &   48 &  48 &  16 &  16 \\
 4 &  4 &  4 &  4 &  8  & &    8  & &   24 &   24 &  24 &  24 &   8 \\
 4 &  4 &  4 &  4 &  4  & &    4  & &   12 &   12 &  12 &  12 &  12 \\
 2 &  4 &  4 &  4 &  4  & &    2  & &   14 &    6 &   6 &   6 &   6 \\
 2 &  2 &  4 &  4 &  4  & &    1  & &    7 &    7 &   3 &   3 &   3 
\end{array}
$$

Then use Lemma \ref{zerostep} with $u_1=u_2=3$, $u_3=u_4=u_5=2$.
\end{Ex}

\begin{Rmk} As in the case of triple Segre products,
in principle, there are no obstructions in repeating the
argument, when we substitute powers of $2$
with powers of $3$ (see the proof of Theorem \ref{main2} we gave in case $a=9$), 
or any other integer $p>1$.

For some numerical cases, the bound for identifiability that we get
using powers of numbers bigger than two, can be closer to
the maximal value $k_{max}$.
\end{Rmk}


\begin{thebibliography}{}
\bibitem[AOP]{AOP}
H.~ABO, G.~OTTAVIANI and C.~PETERSON.
\newblock Induction for secant varieties of {S}egre varieties.
\newblock {\em Trans. Amer. Math. Soc.}, 361(2) (2009) 767--792.

\bibitem[AMR]{AMR}E.~ALLMAN, C. MATIAS and J.~RHODES. 
\newblock Identifiability of parameters in latent structure 
models with many observed variables.
\newblock {\em Ann. Statist},  37  (2009) 3099-3132.

\bibitem[AR]{AR}E.~ALLMAN and J.~RHODES. 
\newblock Phylogenetic invariants for the general Markov model of sequence mutation.
\newblock {\em Mathematical Biosciences}, 186 (2003) 133--144.


\bibitem[tB]{tB} J.~TEN BERGE.
\newblock Partial uniqueness in CANDECOMP/PARAFAC.
\newblock {\em Journal of Chemometrics}, 18 (2004) 12--16.

\bibitem[BCS]{BCS} P. B\"URGISSER, M.CLAUSEN and M.A. SHOKROLLAHI.
\newblock Algebraic Complexity theory.
\newblock{\em Grundl. Math. Wiss.}, 315, Springer, 1997.

\bibitem[CC1]{CC1}
L.~CHIANTINI and C. ~CILIBERTO.
\newblock Weakly defective varieties.
\newblock {\em Trans. Amer. Math. Soc.}, 354(1) (2002) 151--178. 

\bibitem[CC2]{CC2}
L.~CHIANTINI and C. ~CILIBERTO.
\newblock On the $k$-th secant order of a  projective variety.
\newblock {\em J. London Math. Soc.}, 73(2) (2006) 436--454. 

\bibitem[CC3]{CC3}
L.~CHIANTINI and C. ~CILIBERTO.
\newblock On the dimension of secant varieties.
\newblock {\em J. Eur. Math. Soc.}, 12 (2010) 1267--1291.
 
\bibitem[CMO]{CMO}
L.~CHIANTINI, M.~MELLA and G.~OTTAVIANI,
\newblock {\em in preparation}.
 
\bibitem[GS]{GS} D.~GRAYSON and M.~STILLMAN.
\newblock Macaulay 2, a software system for research in algebraic geometry.
\newblock Available at {\tt www.math.uiuc.edu/Macaulay2/ }.

\bibitem[E] {Ein} L.~EIN. 
\newblock Varieties with small dual variety.
\newblock {\em Invent. Math.}, 86 (1986) 63--74.

\bibitem[ERSS]{ERSS} N. ERIKSSON, K. RANESTAD, B. STURMFELS  and 
S. SULLIVANT.
\newblock Phylogenetic Algebraic Geometry. 
\newblock (In C. Ciliberto, A. V. Geramita, B. Harbourne, R. M. Miro-Roig, 
and K. Ranestad (Ed.), Projective Varieties with Unexpected Properties: 
a volume in Memory of Giuseppe Veronese.
\newblock Proceedings of the international conference 
'Varieties with Unexpected Properties', Siena,
Italy, June 8-13, 2004.)
\newblock Walter de Gruyter (2005) 177--197.

\bibitem[KB]{KB}
T.~KOLDA and B.~BADER. 
\newblock Tensor Decompositions and Applications. 
\newblock{\em SIAM Review}, 51(3) (2009) 455--500.

\bibitem[K]{Kru}
J.~B.~KRUSKAL.
\newblock Three-way arrays: rank and uniqueness of trilinear decompositions,
with applications to arithmetic complexity and statistics.
\newblock {\em Lin. Alg. Applic.}, 18(2) (1977) 95--138. 

\bibitem[Land0]{Land0} J.~M.~LANDSBERG. 
\newblock The geometry of tensors with applications.
\newblock {\em in preparation}.

\bibitem[Land1]{Land}
J.~M.~LANDSBERG. 
\newblock Kruskal's theorem.
\newblock {\em Preprint} arXiv:0902.0543v1, 2009, it will appear in chap. 13 of [Land0]. 

\bibitem[Lat]{Lat} L. De LAUTHAWER. 
\newblock A link between the canonical decomposition in multilinear algebra and
simultaneous matrix diagonalization.
\newblock{\em SIAM J. Matrix Anal. Appl.}, 28 (2006) 642--666.

\bibitem[Lick]{Lick} T. LICKTEIG.
\newblock Typical tensorial rank. 
\newblock {\em Linear Algebra Appl.}, 69 (1985) 95--120. 

\bibitem[R]{R}
J.~A.~RHODES. 
\newblock A concise proof of Kruskal's theorem on tensor decomposition.
\newblock {\em Preprint} arXiv:0901.1796, 2009. 

\bibitem[S]{S} 
V. STRASSEN.
\newblock Rank and optimal computation of generic tensors. 
\newblock {\em Linear Algebra Appl.}, 52 (1983) 645--685. 

\bibitem[SB]{SB} N.D. SIDIRIPOULOS and R. BRO. 
\newblock On the uniqueness of multilinear decomposition of N-way arrays. 
\newblock{\em Journal of Chemometrics}, 14, (2000) 229--239.

\bibitem[SS]{St}
A. STEGEMAN and N. D. SIDIROPOULOS.
\newblock On Kruskal's uniqueness for Candecomp/Parafac decomposition. 
\newblock {\em Linear Algebra Appl.}, 420 (2007) 540--552. 

\bibitem[Z]{zak} F. ZAK. 
\newblock Tangents and secants of varieties.
\newblock {\em Transl Math. Monograph}, 127 (1993).

\end{thebibliography}
\end{document}